\documentclass[11pt]{amsart}
\usepackage{epsf,amsfonts,amsmath}
\usepackage{amssymb} 
\usepackage{amsbsy}
\usepackage{amscd} 
\usepackage[dvips]{graphicx}

\newtheorem{theorem}{Theorem}[section]

\theoremstyle{definition}
\newtheorem{definition}[theorem]{Definition}
\newtheorem{example}[theorem]{Example}
\theoremstyle{remark}
\newtheorem{remark}[theorem]{Remark}

\title[Surfaces constructed from groups using
projective planes]{Regular algebraic surfaces isogenous to a higher
  product constructed from group representations using projective
  planes}

\author{Nathan Barker, Nigel Boston, Norbert Peyerimhoff, Alina
  Vdovina} 
\date{$1^{\text{th}}$ August 2011}


\begin{document}

\maketitle

\begin{abstract}
  Regular algebraic surfaces isogenous to a higher product of curves
  can be obtained from finite groups with ramification structures. We
  find unmixed ramification structures for finite groups constructed
  as $p$-quotients of particular infinite groups with special
  presentation related to finite projective planes.  \\ \medskip
\end{abstract}

\section{Introduction}
\label{intro}
 
An algebraic surface is {\em isogenous to a higher product} (of
curves) if it admits a finite unramified covering which is isomorphic
to a product of curves $C_1 \times C_2$ of genera $g(C_i) \ge 2$. It
was shown in \cite{Cat00} that every such surface $S$ has a unique
{\em minimal} realisation $S \cong (C_1 \times C_2)/G$, where $G$ is a
finite group acting freely on $C_1 \times C_2$ and $C_1$ and $C_2$
have the smallest possible genera. Moreover, $G$ respects the product
structure by either acting diagonally on each factor (unmixed case) or
there are elements in $G$ interchanging the factors (mixed
case). Surfaces isogenous to a higher product are always minimal and
of general type. In this paper we restrict our considerations to the
unmixed case.

The {\em irregularity} $q(S)$ of a surface $S$ is the difference
between its geometric and its algebraic genus, and agrees with the
Hodge number $h^{1,0}(S)$. Surfaces with vanishing irregularity are
called {\em regular}. Since $q(S) = g(C_1/G) + g(C_2/G)$ (see
\cite[Prop. 2.2]{Ser}), we have $C_i/G \cong {\mathbb P}^1$ for both
curves in the minimal realisation of a regular surface.

Every surface $S$ isogenous to a higher product gives rise to a finite
group $G$ via its minimal realisation. This process can be reversed.
Starting with a finite group $G$, the existence of a so called {\em
  ramification structure} can be used to construct a regular surface
of the form $(C_1 \times C_2)/G$. We will discuss ramification
structures and the construction of the associated surfaces in Section
\ref{sec:fundnotions}. Bauer, Catanese and Grunewald \cite{BCG} used
this group theoretical description to classify all regular surfaces
$S$ isogenous to a product of curves with vanishing geometric genus
$p_g(S) = h^{2,0}(S)$. The process in \cite{BCG} was aided by the
reduction of the search of ramification structures to groups of order
less than $2000$, for which the MAGMA library of small groups could
then be used. They saw this classification as the solution in a very
special case to the open problem posed by Mumford: \textquotedblleft
Can a computer classify all surfaces of general type with $p_g=0$?''

The infinite group in \cite[Example 6.3]{How} given by the
presentation
\begin{equation} \label{eq:G0}
G_0 := 
\langle x_0,\dots,x_6 \mid x_i x_{i+1} x_{i+3}\, (i \in {\mathbb Z}_7 \rangle
\end{equation}
was used in \cite{BBPV} to construct finite $2$-groups with special
unmixed and mixed ramification structures, giving rise to unmixed and
mixed Beauville surfaces. These finite $2$-groups were the maximal
$2$-quotients of $2$-class $k$ of both the group $G_0$ and its index
two subgroup $H_0$ generated by $x_0$ and $x_1$.

The above group $G_0$ belongs to a family called {\em groups with
  special presentation}. These groups were introduced by Howie
\cite{How} and are related to projective planes over finite fields
(see Section \ref{sec:ourgroups} for more details). It was proved in
\cite{EV} that all groups with special presentation are just infinite
(i.e., they are infinite groups all of whose non-trivial normal
subgroups have finite index). A natural question arose: {\em Do any
  other groups with special presentations give rise to finite groups
  with particular ramification structures?}

In this article we consider finite index subgroups of the groups
listed in \cite[Example 3.3]{EH}, an index $3$ subgroup of the
following group with special presentation from \cite[Example
6.4]{How},
\begin{equation} \label{eq:G}
G := \langle x_0,\dots,x_{12} \mid x_i^3, x_i x_{i+1} x_{i+4}\, 
(i \in {\mathbb Z}_{13}) \rangle
\end{equation}
and the group given in \cite[Example 2]{IV} constructed from a
polyhedral presentation (a generalization of the triangle
presentations defined in \cite{CMSZ1}). We use the computer program
MAGMA to search for unmixed ramification structures in maximal
$p$-quotients of $p$-class $k$ of the above mentioned groups for
various primes $p$.  These ramification structures give then rise to
particular regular surfaces isogenous to a higher product. Our
results are presented in Section \ref{sec:ourgroups} below.

\section{Ramification structures and associated surfaces}
\label{sec:fundnotions}

\subsection{Group theoretical structures}
\label{subsec:ramstruc}

Following \cite{BCG} closely, we give the definition of an (unmixed)
ramification structure of a finite group $G$. 

An $r$-tuple $T = [g_1,\dots,g_r]$ of non-trivial elements of $G$ is
called a {\em spherical system of generators}, if $g_1,\dots,g_r$
generate $G$ and $g_1 g_2 \cdots g_r = 1$. The $r$-tuple
$[m_1,\dots,m_r]$ of non-decreasing orders of the elements $g_i$ is
called the {\em type} of the spherical system $T$ of generators, i.e.,
$2 \le m_1 \le m_2 \le \dots \le m_r$ and there is a permutation $\tau
\in Sym(r)$ such that $m_i = {\rm ord}(g_{\tau(i)})$. Let
$$ \Sigma(T) := \bigcup_{g \in G} \bigcup_{j=0}^\infty \bigcup_{i=1}^r
\{ g g_i^j g^{-1} \} $$ 
be the union of all conjugates of the cyclic subgroups generated by
the elements $g_i$ of the spherical system. Two spherical systems of
generators $T_1=[g_1,\dots,g_r]$ and $T_2=[g_1',\dots,g_s']$ are
called {\em disjoint} if $\Sigma(T_1) \cap \Sigma(T_2) = \{ 1 \}$.  An
unmixed ramification structure is defined as follows.

\begin{definition} (Unmixed ramification structures, see \cite[Definition 1.1]{BCG}) Let $A_1 =
  [m_1,\dots,m_r]$ and $A_2 = [n_1,\dots,n_s]$ be tuples of natural
  numbers with $2 \le m_1 \le \dots \le m_r$ and $2 \le n_1 \le \dots
  \le n_s$. An {\em unmixed ramification structure of type
    $(A_1,A_2)$} for a finite group $G$ is a pair $(T_1,T_2)$ of
  disjoint spherical systems of generators such that $T_1$ has type
  $A_1$ and $T_2$ has type $A_2$.
\end{definition}

The disjointness of the pair $(T_1,T_2)$ of an unmixed ramification
structure guarantees that $G$ acts freely on the product $C_1 \times
C_2$ of associated algebraic curves (see Section
\ref{subsec:groupsurface} and the references therein). In this article
we will only consider unmixed ramification structures and their
associated surfaces. For examples of the mixed case see, e.g.,
\cite{BCG0,BCG,BBPV}.

\subsection{From ramification structures to algebraic surfaces}
\label{subsec:groupsurface}

In this section we explain how to construct an algebraic surface $S =
(C_{T_1} \times C_{T_2})/G$ from a given finite group $G$ with an
unmixed ramification structure $(T_1,T_2)$.

Let $G$ be a finite group and $T=[g_1,\dots,g_r]$ be a spherical
system of generators with $m_i = {\rm ord}(g_{\tau(i)})$. For $1 \le i
\le r$, let $P_1,\dots,P_l \in {\mathbb P}^1$ be a sequence of points
ordered counterclockwise around a base point $P_0$ and $\gamma_i \in
\pi({\mathbb P}^1-\{P_1,\dots,P_r\},P_0)$ be represented by a simple
counterclockwise loop around $P_i$, such that $\gamma_1 \gamma_2 \dots
\gamma_r =1$. By Riemann's existence theorem, we obtain a surjective
homomorphism
$$ \Phi: \pi({\mathbb P}^1-\{P_1,\dots,P_r\},0) \to G $$
with $\Phi(\gamma_i) = g_i$ and a Galois covering $\lambda: C_T \to
{\mathbb P}^1$ with ramification indices equal to the orders of the
elements $g_1,\dots,g_r$. These data induce a well defined action of
$G$ on the curve $C_T$, and by the Riemann-Huritz formula, we have
\begin{equation} \label{eq:gCT}
g(C_T) = 1 + \frac{|G|}{2}\left(r-2-\sum_{l=1}^r \frac{1}{m_l}\right). 
\end{equation}

Now, we assume that $G$ admits an umixed ramification structure
$(T_1,T_2)$. This leads to a diagonal action of $G$ on the product
$C_{T_1} \times C_{T_2}$, and the disjointness of the two spherical
systems of generators ensures that $G$ acts freely on the product of
curves. The associated algebraic surface $S$ is the quotient $(C_{T_1}
\times C_{T_2})/G$. By the Theorem of Zeuthen-Segre, we have for the
topological Euler number
$$ e(S) = 4\frac{(g(C_{T_1})-1)(g(C_{T_2})-1)}{|G|}, $$
as well as the relations (see \cite[Theorem 3.4]{Cat00}),
$$ \chi(S) = \frac{e(S)}{4} = \frac{K_S^2}{8}, $$
where $K_S^2$ is the self intersection number of the canonical divisor
and $\chi(S) = 1 + p_g(S) - q(S)$ is the holomorphic Euler-Poincar\'e
characteristic of $S$. Assume that $(T_1,T_2)$ is of the type
$(A_1,A_2)$ with $A_1=[m_1,\dots,m_r]$ and $A_2=[n_1,\dots,n_s]$. Then
the above relations imply for the associated surface $S$ that
$$ 
\chi(S) = \frac{|G|}{4}\left(r-2-\sum_{l=1}^r \frac{1}{m_l}\right)
\left(s-2-\sum_{l=1}^s \frac{1}{n_l}\right). 
$$

\section{Groups with special presentations}
\label{sec:ourgroups}

As mentioned in \cite{EH}, small cancellation groups are
generalizations of surface groups and satisfy many of the nice
properties of those groups. It was proved in \cite{EH} (with a small
list of exceptions) that almost all groups with a presentation
satisfying the small cancellation conditions $C(3)$ and $T(6)$ contain
a free subgroup of rank $2$. 

Further, \cite{How} proved that most $C(3),T(6)$ groups $G$ (namely,
the ones which do {\em not} have special presentations) are
SQ-universal. (A group $G$ is called SQ-universal if every countable
group can be embedded in a quotient group of $G$.) In that article, a
group presentation was called {\em special} if every relator has
length $3$ and the star graph is isomorphic to the incidence graph of
a finite projective plane. (See \cite[p. 61]{LS} for a textbook
reference on the star graph of a presentation.) Moreover, it was asked
(see \cite[Question 6.11]{How}) whether any or all of the groups with
special presentations are SQ-universal. It was proved in \cite{EV}
that groups with special presentation are just infinite (i.e., all
non-trivial normal subgroups have finite index) and, therefore, cannot
be SQ-universal. (Note that special presentations in the sense of
\cite{How} are $(3,3)$-special in the sense of \cite{EV}.)

Howie \cite{How} also set up an example machine (see Theorem \ref{TSP}
below) to create infinitely many groups with special
presentations. More precisely, he constructed a special presentation
with star graph isomorphic to the incidence graph of the projective
plane over every finite field ${\mathbb F}_q$, where $q$ is a prime
power. Until then, only seven examples of special presentations were
known (see \cite[Example 3.3]{EH}), and each of them has a star graph
isomorphic to the Heawood graph (i.e., the incidence graph of the
$7$-point projective plane over ${\mathbb F}_2$).

Given a finite field $K=\mathbb{F}_q$ (for $q$ a prime power), a
\emph{positive} presentation with star graph isomorphic to this
incidence graph of the Desarguesian projective plane over $K$ is
formed.

The construction takes a cubic extension of $K$, namely
$F=\mathbb{F}_{q^3}$, and identifies the cyclic group
$C_{m}=F^{\times}/K^{\times}$ with the points of the projective pane
$\mathcal{P}$ over $K$, where $m=q^{2}+q+1$.

The group $C_{m}$ acts on $\mathcal{P}$ via multiplication in $F$, and
this action is regular on both the points and lines of $\mathcal{P}$
\emph{i.e.} $C_{m}$ is a Singer group, see \cite{HP}. The lines of
$\mathcal{P}$ can be identified with the subset $\sigma L$ of $C_{m}$,
where $\sigma$ ranges over $C_{m}$ and $L$ is a fixed line or perfect
difference set \emph{i.e.} a set of residues $a_{1},...,a_{q+1}
\mod{m}$ such that every non-zero residue modulo $m=q^{2}+q+1$ can be
expressed uniquely in the form $a_{i}-a_{j}$.
 
\begin{theorem} \cite[Theorem 6.2]{How} \label{TSP}
  Let $q$ be a prime power and $m=q^{2}+q+1$. Then there exists a
  subset $l$ of $q+1$ elements of $\mathbb{Z}_{m}$ such that
  $$ 
  \langle x_{0},...,x_{m-1} \mid x_{i}x_{i+\lambda}x_{i+\lambda+q\lambda} \  
  (i \in \mathbb{Z}_m, \lambda \in l) \rangle
  $$
  is a special presentation whose star graph is isomorphic to the
  incidence graph of the projective plane over $GF(q)$.
\end{theorem}

Let us now present results on ramification structures of finite groups
obtained from particular groups $G$ with special
representations. These finite groups are generated via the lower,
exponent $p$-central series, i.e.,
$$ G = P_0(G) \ge \dots \ge P_{i-1}(G) \ge P_i(G) \ge \dots, $$
where $P_i(G) = [P_{i-1}(G),G]P_{i-1}(G)^p$ for $i \ge 1$. The finite
groups $G_{p,k}$ under considerations are then the maximal
$p$-quotients of $p$-class $k$, denoted by $G_{p,k}$ and given by
$G_{p,k} = G/P_k(G)$.

The results discussed below are obtained via the computer program
MAGMA (see \cite{magma}). Note that the algorithm \texttt{pQuotient}
constructs, for a given group $G$, a consistent power-conjugate
presentation for $G_{p,k}$.

\subsection{Ramification structures for the group in \cite[Example
  6.3]{How}}
\label{subsec:ex6.3How}

The group $G_0$ in \eqref{eq:G0} with seven generators $x_0,\dots,x_6$
appeared as Example 6.3 in \cite{How}. ($G_0$ is constructed using
Theorem \ref{TSP} with $q=2$, $m=7$ and $l=\{1,2,4\}$.) The subgroup
$H_0$ generated by $x_0, x_1$ has index two. In \cite[Theorems 4.1 and
4.2]{BBPV}, we presented unmixed ramification structures for the
$2$-groups $(H_0)_{2,k}$ for $3 \le k \le 64$ (for $k$ not a power of
$2$), as well as mixed ramification structures for the $2$-groups
$(G_0)_{2,k}$ for $3 \le k \le 10$ (again, for $k$ not a power of
$2$). Since the involved spherical systems of generators consist of
three elements, these ramification structures are actually Beauville
structures and lead to new examples of Beauville surfaces.

The group $G_0$ appears also in \cite[Section 4]{CMSZ2} as the group
A.1. (The articles \cite{CMSZ1,CMSZ2} are concerned with simply
transitivie group actions on the vertices of ${\widetilde
  A}_2$-buildings.) The index two subgroup $H_0$ was also used in
\cite{PV} to construct families of expander graphs of vertex degree
four.

\subsection{Ramification structures for the groups in \cite[Example
  3.3]{EH}}
\label{subsec:ex3.3EH}

There, a list of seven special group presentations $G_i = \langle {\bf
  x} \mid {\bf r}_i \rangle$, $1 \le i \le 7$ were given
($\mathbf{x}=\{a,b,c,d,e,f,g\}$). The star graphs of all seven
presentations are isomorphic to the incidence graph of the $7$-point
projective plane. Our group $G_0$ in \eqref{eq:G0} coincides with
their group $G_3$, which was discussed in Section
\ref{subsec:ex6.3How}. It is stated in \cite{EH} that the only
isomorphism between abelianised groups $G_i^{ab}$ is between
$G_4^{ab}$ and $G_6^{ab}$. However, if one looks at the commutator
subgroup $C_{4}$ and $C_{6}$ of the groups $G_{4}$ and $G_{6}$, then
$C_{4}^{ab}\cong \mathbb{Z}/4\mathbb{Z}$ and $C_{6}^{ab}\cong
\mathbb{Z}/2\mathbb{Z}$. Therefore, $G_{4}$ can not be isomorphic to
$G_{6}$.

We use the computer program MAGMA to search in the maximal
$p$-quotients of maximal class $k$ for $1\le k\le 10$ of certain
finite index subgroups of the groups $G_{i}$ for unmixed ramification
structures.

\subsubsection{Special presentation $G_{1}$}
There is a subgroup $H_{1}$ of index $4$ in $G_{1}$ generated by
$b$. Thus, as $H_{1}\cong \mathbb{Z}$ all maximal $p$-quotients of
$p$-class $k$ of $H_{1}$ are cyclic groups of order $p^k$. Therefore,
there will be no unmixed ramification structures coming from the
groups $(H_{1})_{p,k}$.

\subsubsection{Special presentation $G_{2}$ and $G_{4}$}
There is a subgroup $H_{2}$ of index $16$ in $G_{2}$ (the commutator
subgroup of $G_{2}$) generated by $h_{2,0}=bd^{-1}a^{-1}bc^{-1}$,
$h_{2,1}=abd^{-1}abd^{-1}bc^{-1}$, $h_{2,2}=(bc^{-1})^{2}$,
$h_{2,3}=adc^{-1}$ and $h_{2,4}=db^{-1}a^{-1}$
$db^{-1}a^{-1}bc^{-1}$. The maximal $7$-quotients of $7$-class $k$,
written as $(H_{2})_{7,k}$, gives rise to a disjoint pair of spherical
generators of length $3$ given by the tuples,
$$T_{2,1}=[h_{2,0},h_{2,1},h_{2,1}^{-1}h_{2,0}^{-1}] \ \text{and} \ 
T_{2,2}=[h_{2,0}h_{2,1}^{2},h_{2,0}h_{2,1}^{3},
h_{2,1}^{-3}h_{2,0}^{-1}h_{2,2}^{-2}h_{2,0}^{-1}],$$ 
for $1\le k\le 10$. To simplify notation, we denoted the images of
$h_{2,0}$ and $h_{2,1}$ in $(H_{2})_{7,k}$, again, by $h_{2,0}$ and
$h_{2,1}$. Thus, the groups $(H_{2})_{7,k}$ have unmixed ramification
structures.

The group $G_4$ of \cite[Example 3.3]{EH} coincides with the group
C.1. in \cite[Section 5]{CMSZ2} (via the identification $a_0=a, a_1=f,
a_2=c, a_3=d, a_4=e, a_5=g, a_6=b$). We find a subgroup $H_{4}$ of
index 48 in $G_{4}$ generated by $h_{4,0}=da^{-1}bc^{-1}$,
$h_{4,1}=bc^{-1}bc^{-1}ea^{-1}bf^{-1}$,
$h_{4,2}=cb^{-1}ae^{-1}cf^{-1}$, $h_{4,3}=(ea^{-1}bf^{-1})^{2}$,
$h_{4,4}=cb^{-1}cb^{-1}ea^{-1}bf^{-1}$ and
$h_{4,5}=ea^{-1}bc^{-1}ea^{-1}bc^{-1}ea^{-1}bf^{-1}$. For $1\le k\le
10$, the maximal $7$-quotients of $7$-class $k$, $(H_{4})_{7,k}$,
gives rise to a disjoint pair of spherical generators of length $3$
given by the tuples,
$$T_{4,1}=[h_{4,0},h_{4,1},h_{4,1}^{-1}h_{4,0}^{-1}] \ \text{and} \ 
T_{4,2}=[h_{4,0}h_{4,1}^{2},h_{4,0}h_{4,1}^{3},
h_{4,1}^{-3}h_{4,0}^{-1}h_{4,2}^{-2}h_{4,0}^{-1}].$$ The groups
$(H_{2})_{7,k}$ and $(H_{4})_{7,k}$ have the same maximal
$7$-quotients of maximal $7$-class $k$ for $1\le k\le 10$. However,
the abelianizations of the infinite groups are $H_{2}^{ab}\cong
\mathbb{Z}_{7}\times \mathbb{Z}_{21}$ and $H_{4}^{ab}\cong
\mathbb{Z}_{7}\times \mathbb{Z}_{7}$. The following theorem summarizes
the above observations.

\begin{theorem} 
  For $r=2,4$, $k=1,\dots,10$, the groups $(H_{r})_{7,k}$ are of order
  $7^{a}$ and admit unmixed ramification structures
  $(T_{r,1},T_{r,2})$ of type $([7^{b},7^{b},7^{b}],$
  $[7^{b},7^{b},7^{b}])$ for
    $$ a= \begin{cases} 2k & \text{if $k= 1,2,5,8,9$}, \\ 2k-1 &
      \text{if $k=2,3,6,7,10.$} \end{cases}\ \text{and} \
    b= \begin{cases} 2 & \text{if $1\le k\le 4$}, \\ 3 & \text{if
        $5\le k\le 8$}, \\ 5 & \text{if $9\le k\le 10$}. \end{cases}$$
\end{theorem}

The unmixed ramification structures given for the groups
$(H_{2})_{7,k},(H_{4})_{7,k}$ above give rise to unmixed Beauville
surfaces $S=(C_{T_1} \times C_{T_2})/(H_{n})_{7,k}$ for $n=2,4$. For
example, the order of the group $(H_{2})_{7,1}$ and $(H_{4})_{7,1}$ is
$7^2$. Therefore, the genera of the curves $C_{T_i}$ is (see
\eqref{eq:gCT})
$$ g(C_{T_1}) = g(C_{T_2}) = 1 + 2\times7 =15 , $$
and the holomorphic Euler-Poincar{\'e} characteristic of $S$ is
$$ \chi(S) = \frac{(g(C_{T_1})-1)(g(C_{T_2})-1)}{|G|} = 4. $$

\begin{remark}(see \cite[Beauville's examples 3.22]{Cat00}) 
  The groups $(H_{2})_{7,1},(H_{4})_{7,1}$ are isomorphic to the group
  $(\mathbb{Z}/7\mathbb{Z})^{2}$ and the two curves
  $C_{T_{1}}=C_{T_{2}}$ are given by the Fermat curve $x^7+y^7+z^7=0$
  of degree $7$. The group $(\mathbb{Z}/7\mathbb{Z})^{2}$ acts on
  $C_{T_{1}}\times C_{T_{2}}$ by the following rule
  \[(\alpha,\beta) \cdot
  ([x:y:z],[u:v:w])=([\xi^{\alpha}x:\xi^{\beta}y:z],
  [\xi^{\alpha+2\beta}u:\xi^{\alpha+3\beta}v:w]),\] 
  where $\xi=e^{\frac{2\pi i}{7}}$ and $\alpha,\beta \in
  \mathbb{Z}/7\mathbb{Z}$. We identify $h_{n,0}\mapsto \alpha$ and
  $h_{n,1}\mapsto \beta$ for $n=2,4$.
\end{remark}

\subsubsection{Special presentation $G_{5}$}
The group $G_{5}$ coincides with the group A.2 in \cite[Section
5]{CMSZ2}. We find a subgroup $H_{5}$ of index $3$, generated by
$h_{5,0}=ba^{-1}$, $h_{5,1}=ca^{-1}$, $h_{5,2}=da^{-1}$,
$h_{5,3}=ea^{-1}$, $h_{5,4}=fa^{-1}$ and $h_{5,5}=ga^{-1}$ which have
the same maximal $2$-quotients of $2$-class $k$ as the group $G_{0}$
in (\ref{eq:G0}) for $1\le k \le 10$. However, the abelianization of
this group is $H_{5}^{ab}\cong
\mathbb{Z}_{2}\times\mathbb{Z}_{2}\times\mathbb{Z}_{14}$ which is not
isomorphic $G^{ab}_{0}\cong
\mathbb{Z}_{2}\times\mathbb{Z}_{2}\times\mathbb{Z}_{6}$.

In addition, we have a subgroup $F_{5}$ in $H_{5}$ of index $2$, which
appears to have the same maximal $2$-quotients of $2$-class $k$ as
$H_{0}$ (the subgroup of $G_{0}$ generated by $x_{0},x_{1}$) for $1\le
k \le 10$. The abelianization of this group is $F_{5}^{ab}\cong
\mathbb{Z}_{4}\times\mathbb{Z}_{28}$ which is not isomorphic
$G^{ab}_{0}\cong \mathbb{Z}_{4}\times\mathbb{Z}_{12}$.

\subsubsection{Special presentation $G_{6}$ and $G_{7}$}
The group $G_{6}$ coincides with the group B.2 in \cite[Section
5]{CMSZ2}. We have the group specified by relations $\mathbf{r}_{6}=$
$\{abe,acb,aec,bf^2,cd^2,dfg,eg^2\}$ on $7$ generators but can be
rewritten to a group generated by $\{a,b\}$ with relations
$\mathbf{r}^{'}_{6}=\{b^{-1}a^{-1}b^{2}a^{-2}b^{-3}a^{-1},a^{3}baba^{-2}b^{2}\}$
(see \cite[Section 2.7]{SBT}). The group $G_{7}$ coincides with the
group B.1 in \cite[Section 5]{CMSZ2}.

We see that both groups have a subgroup $H_{6},H_{7}$ of index 24 in
$G_{6},G_{7}$, respectively, which gives rise to maximal $3$-quotients
of $3$-class $k$ for $1\le k\le 10$. However, the $3$-groups are too
large to successfully search for unmixed ramification structures. The
abelianization of both groups is $H_{6}^{ab}\cong H_{7}^{ab}\cong
\mathbb{Z}_{3}\times\mathbb{Z}_{3}\times\mathbb{Z}_{3}\times\mathbb{Z}_{3}$.

\subsection{Ramification structures for the group in \cite[Example
  6.4]{How}}
\label{subsec:ex6.4How}

The group $G$ in \eqref{eq:G} with $13$ generators $x_0,\dots,x_{12}$
appeared as Example 6.4 in \cite{How}. ($G$ is constructed using
Theorem \ref{TSP} with $q=3$, $m=13$ and $l=\{0,1,3,9\}$.) The
subgroup $H$ generated by $x_0, x_1, x_2$ has index $3$. Again, the
group $G$ can also be found in \cite[Section 4]{CMSZ2} as the group
1.1 (via the identification $a_i = x_{2i}$, where the indices are
taken modulo $13$).

For simplicity, let the elements in the finite $3$-quotients $H_{3,k}$
corresponding to $x_0,x_1,x_2 \in H$ be denoted, again, by
$x_0,x_1,x_2$. Let $x = x_2^{-1} x_1^{-1} x_0^{-1}$, $y_0=x_0 x_1^2
x_2^2$, $y_1=x_0^2 x_1 x_2^2$, $y_2=x_1 x_2^{-1} x_2^{x_0}$ and $y =
y_2^{-1} y_1^{-1} y_0^{-1}$. Using MAGMA, we obtain the following
result.

\begin{theorem} \label{thm:main}
  For $k=2,\dots,60$, the groups $H_{3,k}$ are of order $3^{a_k}$ and
  admit unmixed ramification structures $(T_1,T_2)$ of type
  $([3,3,3,3^{d_k}],[3^{b_k},3^{b_k},3^{b_k},3^{b_k}])$, where $T_1 =
  (x_0,x_1,x_2,x)$, $T_2 = (y_0,y_1,y_2,y)$, $b_k = 1+[\log_3
  \frac{3k}{4}]$, $d_k = 1+[\log_3 k]$, and
  $$ a_k = \begin{cases} 8j & \text{if $k =3j$}, \\ 8j+3 &
    \text{if $k=3j+1$}, \\ 8j+6 & \text{if $k=3j+2$}. \end{cases}$$
  Here $[x]$ denotes the largest integer $\le x$.
\end{theorem}

We conjecture that this result holds true for all integers $k \ge 2$.

\medskip

The ramification structures of $H_{3,k}$ in Theorem \ref{thm:main}
give rise to algebraic surfaces $S=(C_{T_1} \times
C_{T_2})/H_{3,k}$. For example, the order of the group $H_{3,2}$ is
$a_2=3^6$. Therefore, the genera of the curves $C_{T_i}$ is (see
\eqref{eq:gCT})
$$ g(C_{T_1}) = g(C_{T_2}) = 1 + 3^5 = 244, $$
and the holomorphic Euler-Poincar{\'e} characteristic of $S$ is
$$ \chi(S) = \frac{(g(C_{T_1})-1)(g(C_{T_2})-1)}{|G|} = 81. $$

\subsection{Ramification structures for the groups of Theorem
  \ref{TSP} with $q\ge 4$}
The construction given by Theorem \ref{TSP} is for any $q$ a prime
power. For $q=4$ the group below is given.

\begin{example}\cite[Example 6.5]{How}
  We have that $q^{2}+q+1=21$ and so
  $\mathbb{F}^{\times}_{q^{3}}/\mathbb{F}^{\times}_{q}$ is identified
  with $\mathbb{Z}_{21}$. The group $\widehat G$ is given by the
  presentation,
  \begin{equation}\label{groupq=4}
    {\widehat G}:=\langle x_{0},...,x_{20}|
    x_{i}x_{i+7}x_{i+14},x_{i}x_{i+14}x_{i+7},x_{i}x_{i+3}x_{i+15} \ 
    \text{for} \ i\in\mathbb{Z}_{21} \rangle.
  \end{equation}
\end{example}

The abelianization of this group is ${\widehat G}^{ab}\cong
\mathbb{Z}_{2}\times \mathbb{Z}_{2}\times \mathbb{Z}_{2}\times
\mathbb{Z}_{2}\times\mathbb{Z}_{6}\times\mathbb{Z}_{6}$. The group has
maximal $2$-quotients of $2$-class $k$ for $1\le k\le 10$. However, it
is extremely difficult to search for ramification structures of the
maximal $p$-quotients of $p$-class $k$ for $q\ge 5$. The finite groups
are too large and have too many conjugacy classes, which leads to a
computational expensive search.

\subsection{Ramification structures for the group in \cite[Example
  2]{IV}}

In \cite{IV}, a new construction of groups presentations based on
finite projective planes was introduced, generalizing the triangle
presentations of \cite{CMSZ1,CMSZ2}. For the readers convenience, we
explain this briefly. The construction is based on the following
general definition.

\begin{definition} (see \cite{IV}) Let
  $\mathcal{P}_{1},...,\mathcal{P}_{n}$ be $n$ disjoint finite
  projective planes of order $q$. Let $P_{i}$ and $L_{i}$ be the sets
  of points and lines respectively in $\mathcal{P}_{i}$. Let $P=\cup
  P_{i}$, $L=\cup L_{i}$, $P_{i}\cap P_{j}=\emptyset$ for $i\ne j$ and
  let $\lambda$ be a bijection $\lambda : P\to L$.

  A set $\mathcal{K}$ of $k$-tuples $(x_{1},...,x_{k})$ will be called
  a {\em polyhedral presentation} over $P$ compatible with $\lambda$ if
  \begin{enumerate}
  \item given $x_{1},x_{2}\in P$ then $(x_{1},...,x_{k})\in
    \mathcal{K}$ for some $x_{3},...,x_{k}$ if and only if $x_{2}$ and
    $\lambda(x_{1})$ are incident;
  \item $(x_{1},...,x_{k})\in \mathcal{K}$ implies that
    $(x_{2},...,x_{k},x_{1})\in \mathcal{K}$;
  \item given $x_{1},x_{2}\in P$, then $(x_{1},...,x_{k})\in
    \mathcal{K}$ for at most one $x_{3}\in P$.
  \end{enumerate}
  We call $\lambda$ a {\em basic bijection}.
\end{definition}

A polyhedral presentation $\mathcal{K}$ gives rise to a group
presentation $G_{\mathcal{K}}$ in the following way: the generators of
$G_{\mathcal{K}}$ are given by $\cup P_{i}$ and the relations are the
$k$-tuples of $\mathcal{K}$, each written as a product.

\begin{example}
  The {\em triangle presentations} listed in \cite{CMSZ2} can be seen as
  special cases of polyhedral presentations for $n=1$,
  $k=3$ and $q=2,3$.

  We now discuss the case $n=1$, $q=2$. We enumerate the points of the
  projective plane by $1,2,\dots,6$. The following array illustrates a
  basic bijection $\lambda$:
  \[\begin{matrix}
    0: & 1 & 4 & 2\\
    1: & 3 & 2 & 5\\
    2: & 4 & 3 & 6\\
    3: & 0 & 4 & 5\\
    4: & 1 & 5 & 6\\
    5: & 0 & 2 & 6\\
    6: & 0 & 1 & 3.
  \end{matrix}\] 
  Here, every point $k$ represents a row and is followed by the points
  contained in the associated line $\lambda(k)$. For example, the line
  $\lambda(3)$ consists of the points $0,4,5$.

  A triangle presentation $\mathcal{T}$ for the group A.1 in
  \cite{CMSZ2} is given by
  \[(0,1,3),(1,2,4),(2,3,5),(3,4,6),(4,5,0),(5,6,1),(6,0,2),\] and all
  the cyclic permutations, \emph{i.e.} for $(0,1,3) \in \mathcal{T}$
  we also have $(1,3,0)$, $(3,0,1) \in \mathcal{T}$. The associated
  group presentation $G_{\mathcal T}$ agrees with the presentation of
  $G_0$ in \eqref{eq:G0}.
\end{example}

\begin{example}\cite[Example 2]{IV}
  The projective plane $\mathcal{P}$ of order $4$ can be partitioned
  by three projective planes of order two (see \cite{Bruck}). We
  denote points of the subplane $\mathcal{P}_{i}$ for $i=1,2,3$ by
  numbers from $7i-6$ to $7i$. Note that lines in $\mathcal{P}$
  consist of five points, while the lines in $\mathcal{P}_{i}$ consist
  of three points. A basic bijection $\lambda$ for $\mathcal{P}$ is
  given below. Note that each subplane $\mathcal{P}_{i}$ has its own
  basic bijection, denoted by $\lambda_{i}$, satisfying
  $\lambda_{i}(k) \subset \lambda(k)$. In the array below, the row
  associated to the point $k$ lists first the three points in the
  associated line via the basic bijection in the subplane, followed up
  by the two remaining points in $\lambda(k)$.

\begin{minipage}[b]{0.4\linewidth}
\[\begin{matrix}
4: & 5 & 6 & 7 & & 12 & 18\\
7: & 1 & 2 & 5 & & 8 & 21\\
2: & 3 & 4 & 5 & & 14 & 16\\
5: & 1 & 3 & 6 & & 10 & 19\\
1: & 2 & 4 & 6 & & 9 & 15\\
3: & 1 & 4 & 7 & & 11 & 17\\
6: & 2 & 3 & 7 & & 13 & 20
\end{matrix}\]
\end{minipage}
\begin{minipage}[b]{0.2\linewidth}
\phantom{Lala}
\end{minipage}
\begin{minipage}[b]{0.4\linewidth}
\[\begin{matrix}
9: & 12 & 13 & 14 & & 1 & 15\\
11: & 8 & 9 & 12 & & 3 & 17\\
14: & 10 & 11 & 12 & & 2 & 16\\
12: & 8 & 10 & 13 & & 4 & 18\\
10: & 9 & 11 & 13 & & 5 & 19\\
13: & 8 & 11 & 13 & & 6 & 20\\
8: & 9 & 10 & 14 & & 7 & 21
\end{matrix}\]
\end{minipage}

\[\begin{matrix}
18: & 19 & 20 & 21 & & 4 & 12\\
21: & 15 & 16 & 19 & & 7 & 8\\
16: & 17 & 18 & 19 & & 2 & 14\\
19: & 15 & 17 & 20 & & 5 & 10 \\
15: & 16 & 18 & 20 & & 1 & 9\\
17: & 15 & 18 & 21 & & 3 & 11\\
20: & 16 & 17 & 21 & & 6 & 13
\end{matrix}\]

The above basic bijections give rise to the following polyhedral
presentation $\mathcal{K}$ for a projective plane of order $4$,
induced by polyhedral presentations of projective planes of order $2$
\begin{multline*}
(1,9,15),(1,15,9),(2,14,16),(2,16,14),(3,11,17),(3,17,11),(4,12,18),\\
(4,18,12),(5,10,19),(5,19,10),(6,13,20),(6,20,13),(7,8,21),(7,21,8),\\
(1,2,3),(1,4,5),(1,6,7),(3,4,6),(3,7,5),(2,5,6),(2,4,7),(8,9,12),(8,10,13),\\
(8,14,11),(9,14,10),(9,13,11),(12,13,14),(10,11,12),(15,16,17),(15,18,19),\\
\hfill (17,18,20),(17,21,19),(16,19,20),(16,18,21), \hfill
\end{multline*}
and all their cyclic permutations.

All relators in the group presentation $G_{\mathcal{K}}$
given by $\mathcal{K}$ are of length $3$ and the star graph is
isomorphic to the incidence graph of a finite projective plane of
order $4$. This means, by \cite{How}, that the group given by this
presentation $G_{\mathcal{K}}$ is a {\em special presentation}. It
also can be seen that this group acts on a Euclidean building where
the vertex links are incidence graphs of projective planes of order
$4$, see \cite{IV}.

\end{example}

\medskip

There are remarkable connections between the group $G_{\mathcal{K}}$
and the group $G_{0}$ discussed in Subsection
\ref{subsec:ex6.3How}. Firstly, we can present the group
$G_{\mathcal{K}}$ in an alternative way with different generators:
\begin{multline}
  G_{\mathcal K} =  \langle w_{0},...,w_{6},y_{0},...,y_{6},z_{0},...,z_{6} \ | \  
  w_{i}w_{i+1}w_{i+3},y_{i}y_{i+1}y_{i+3},z_{i}z_{i+1}z_{i+3},\\
  \hfill w_{i}^{-1}y_{6(1+i)}z_{i}^{-1},w^{-1}_{i}z_{i}^{-1}y_{6(1+i)}
  (i \in \mathbb{Z}_{7})\rangle, \hfill \label{altdescr}
\end{multline}
where each of the three subsets of generators has very similar
relators like those appearing for the group $G_0$ in \eqref{eq:G0},
with only two more series of relators added representing connections
between the generators of different subsets.

Secondly, the maximal $2$-quotients of $2$-class $k$ of the group
$G_{\mathcal K}$ are isomorphic to the maximal $2$-quotients of
$2$-class $k$ for the group $G_{0}$ (given by the presentation in
\eqref{eq:G0}) for $1\le k\le 20$. However, the groups $G_{\mathcal
  K}$ and $G_{0}$ are not isomorphic, as they have different
abelianized groups $G_{\mathcal K}^{ab}\cong \mathbb{Z}_{2}\times
\mathbb{Z}_{6}\times\mathbb{Z}_{6}$, while $G_{0}^{ab}\cong
\mathbb{Z}_{2}\times \mathbb{Z}_{2}\times\mathbb{Z}_{6}$.
  
\begin{remark}
  If we replace the relators
  $w_{i}^{-1}y_{6(1+i)}z_{i}^{-1},w^{-1}_{i}z_{i}^{-1}y_{6(1+i)} $ in
  \eqref{altdescr} by the relators $x_{i}y_{i}z_{i}$,
  $x_{i}z_{i}y_{i}$ we obtain a group $G^{'}$ with the following presentation
\begin{multline*}
  G^{'} = \langle w_{0},...,w_{6},y_{0},...,y_{6},z_{0},...,z_{6} \ |
  \ w_{i}w_{i+1}w_{i+3},y_{i}y_{i+1}y_{i+3},z_{i}z_{i+1}z_{i+3},\\
  \hfill w_{i}y_{i}z_{i},w_{i}z_{i}y_{i} (i \in
  \mathbb{Z}_{7})\rangle. \hfill
\end{multline*}
This group $G^{'}$ is isomorphic to the group $\widehat G$ given by
the presentation (\ref{groupq=4}) under the identifications

\begin{minipage}[b]{0.3\linewidth}
\[\begin{matrix}
x_{0} \mapsto&z_{0}^{-1}, \\
x_1 \mapsto&w_1^{-1}, \\
x_2 \mapsto&y_2^{-1}, \\
x_3\mapsto &z_3^{-1}, \\
x_4\mapsto&w_4^{-1},\\
 x_5 \mapsto&y_5^{-1}, \\
 x_6 \mapsto&z_6^{-1},
\end{matrix}\]
\end{minipage}
\begin{minipage}[b]{0.3\linewidth}
\[\begin{matrix}
 x_7 \mapsto &w_0^{-1},\\
  x_8 \mapsto&y_1^{-1}, \\
  x_9 \mapsto&z_2^{-1}, \\
  x_{10} \mapsto &w_3^{-1},\\
   x_{11} \mapsto &y_4^{-1}, \\
   x_{12} \mapsto &z_5^{-1}, \\
   x_{13} \mapsto &w_6^{-1},
 \end{matrix}\]
\end{minipage}
 \begin{minipage}[b]{0.3\linewidth}
\[\begin{matrix}
  x_{14} \mapsto &y_0^{-1}, \\
  x_{15} \mapsto &z_1^{-1}, \\
  x_{16} \mapsto &w_2^{-1}, \\
  x_{17} \mapsto &y_3^{-1}, \\
  x_{18} \mapsto &z_4^{-1}, \\
  x_{19}\mapsto &w_5^{-1}, \\
  x_{20} \mapsto &y_6^{-1}.
  \end{matrix}\]
\end{minipage}
\end{remark}

\section*{Appendix: A representation for the group G}\label{Rep}

We include a representation for the group $G$ (given by \eqref{eq:G})
in ${\rm GL}(9,\mathbb{F}_{3}[1/Y])$, which may be useful in the
future (as the matrix representations for the group $G_{0}$ with
presentation (\ref{eq:G0}) were useful for several works
\cite{CMSZ1,PV,LSV}). The representation is due to Donald Cartwright
and the algebra program REDUCE. Recall that the group $G$ coincides
with the group $1.1$ in \cite{CMSZ2}, where we relate the generators
by $a_{i}=x_{2i}$ for $i=0,...,12$, with indices taken modulo $13$. We
set

{\Small
\[
x_{0}: \left(\begin{matrix}
      1 & 1 & 1 & 0 & 2& 2&0 &1 &1 \\
      0 & 1 & 2 & 0 & 0& 1&0 &0 &2 \\       
      0 & 0 & 1 & 0 & 0& 0&0 &0 &0 \\ 
      0 & 0 & 0 & 1 & 1& 1&0 &2 &2 \\
      0 & 0 & 0 & 0 & 1& 2&0 &0 &1 \\       
      0& 0 & 0 & 0 & 0& 1&0 &0 &0 \\ 
       0 & 0 & 0 & 0 & 0& 0&1 &1 &1\\
      0 & 0 & 0 & 0 & 0& 0&0 &1 &2 \\       
      0 & 0 & 0 & 0 & 0& 0&0 &0&1 \\ 
         \end{matrix}\right) +\frac{1}{Y}
\left(\begin{matrix}
     0 & 2 & 2 & 0 & 1& 1&0 &2 &2 \\
      0 & 0 & 1 & 0 & 0& 2&0 &0 &1\\       
      0 & 0 & 0 & 0 & 0& 0&0&0 &0 \\ 
       0 & 1 & 1 & 0 & 2& 2&0 &1 &1 \\
      0 & 0 & 2 & 0 & 0& 1&0 &0 &2 \\       
      0 & 0 & 0 & 0 & 0& 0&0 &0 &0 \\ 
       0 & 2 & 2 & 0 & 1& 1&0 &2 &2 \\
      0 & 0 & 1 & 0 & 0& 2&0 &0 &1 \\       
      0 & 0 & 0 & 0 & 0& 0&0 &0 &0 \\  
         \end{matrix}\right)\]
and 
\[\tau:
\left(\begin{matrix}
      1 & 0 & 0 & 0 & 0& 0&0 &0 &0 \\
      0 & 1 & 0 & 0 & 0& 0&0 &0 &0 \\       
      0 & 0 & 1 & 0 & 0& 0&0 &0 &0 \\ 
      0 & 0 & 0 & 0 & 1& 0&0 &0 &0 \\
      0 & 0 & 0 & 0 & 1& 1&0 &0 &0 \\       
      0 & 0 & 0 & 1 & 0& 1&0 &0 &0 \\ 
      0 & 0 & 0 & 0 & 0& 0&2 &2 &0 \\
      0 & 0 & 0 & 0 & 0& 0&0 &1 &2 \\       
      0 & 0 & 0 & 0 & 0& 0&2 &0 &1 \\ 
         \end{matrix}\right),\] 
}

\noindent
where the other generators $x_{1},...,x_{12}$ are formed via
conjugation of $x_{0}$ by $\tau$, \emph{i.e.}
$x_{i}=\tau^{i}x_{0}\tau^{-i}$ for $i=1,...,12$.

The idea in creating this representation is to write
$\mathbb{F}_{27}=\mathbb{F}_{3}(\theta)$, where $\theta$ is a
primitive element on $\mathbb{F}_{27}$ satisfying
$\theta^{3}=\theta+1$, and to use the basis
$\{\theta^{i}\sigma^{j}|i,j=0,1,2\}$ for the divison algebra
$\mathcal{A}$ over $\mathbb{F}_{27}(Y)$ for an indeterminate $Y$ (in
the order $1,\theta, \theta^{2}, \sigma, \theta\sigma,
...,\theta^2\sigma^2$). Here $\sigma$ is assumed to satisfy $\sigma^3
= Y - 1$ (which implies $(1+\sigma)^{-1}=(1/Y)(1-\sigma+\sigma^{2})$)
and $\sigma \theta \sigma^{-1} = \theta^3$. The generators of
$\mathcal{T}_{\mathcal{K}}$, where $\mathcal{K}$ is a triangle
presentation from \cite{CMSZ1,CMSZ2}, are the
$a_{u}=u^{-1}(1+\sigma)u$, where $u\in
\mathbb{F}_{27}^{\times}/\mathbb{F}_3^{\times}$. Since
$\mathbb{F}_{27}^{\times} = \mathbb{F}_3^{\times}\cdot \{
1=\theta^{13},\theta,\dots,\theta^{12} \}$, we choose
$\alpha_k=\theta^{-k}(1+\sigma)\theta^k$ as in
\cite[p. 178]{CMSZ2}. The $\alpha_k$'s act on $\mathcal{A}$ by
conjugation. A straightforward calculation yields
\begin{multline*}
\alpha_k \theta^i\sigma^j \alpha_k^{-1} =
\theta^i \sigma^j \frac{1}{Y} + \left( \theta^{3i+2k} - \theta^{i+2
    \cdot 3^j k} \right)\sigma^{j+1} \frac{1}{Y} \\ + \left( \theta^{i+8
    \cdot 3^j k} - \theta^{3i+2k+2\cdot 3^{j+1} k} \right) \sigma^{j+2}
\frac{1}{Y} + \theta^{3i+2k+8 \cdot 3^{j+1}k} \sigma^j
\frac{Y-1}{Y}.
\end{multline*} 
Expressing the conjugation by $\alpha_k$ with respect to the above
basis of $\mathcal{A}$ then gives rise to a representation as a $9
\times 9$ matrix over the field $\mathbb{F}_{3}(1/Y)$. We conclude
from \cite{CMSZ2} that the matrices associated to the $\alpha_k$
satisfy the relations of our generators $x_k$. Note, finally, that the
above matrix for $\tau$ represents the conjugation by $\theta$ in
$\mathcal{A}$, i.e., $z \mapsto \theta^{-1} z \theta$.

\bigskip

{\bf Acknowledgements:} We thank Donald Cartwright for the
representations and method given in the Appendix and helpful
correspondences. The first author also wishes to thank Uzi Vishne for
useful correspondences. The research of Nigel Boston is supported by NSA Grant MSN115460.

\medskip

Nathan Barker and Alina Vdovina\\
School of Mathematics and Statistics\\
Newcastle University,\\
Newcastle upon Tyne, NE1 7RU, UK\\
email: nathan.barker@ncl.ac.uk, alina.vdovina@ncl.ac.uk

\medskip

Nigel Boston\\
Department of Mathematics,\\ 
303 Van Vleck Hall, \\
480 Lincoln Drive, \\
Madison, WI 53706, USA\\
email: boston@math.wisc.edu

\medskip

Norbert Peyerimhoff\\
Department of Mathematical Sciences\\
Durham University\\
Science Laboratories\\
South Road, Durham, DH1 3LE, UK\\
email: norbert.peyerimhoff@durham.ac.uk

\medskip


\begin{thebibliography}{99}

\bibitem{BBPV} N. Barker, N. Boston, N. Peyerimhoff \and A. Vdovina,
  New examples of Beauville surfaces, to appear in {\em
    Monatsh. Math.}  (2011).

\bibitem{SBT} S. Barre, Poly\'{e}dres de rang deux, Thesis ENS Lyon, December (1996) http://web.univ-ubs.fr/lmam/barre/these1.pdf.

\bibitem{BCG0} I.C. Bauer, F. Catanese \and F. Grunewald, Beauville
  surfaces without real structures, in {\em Geometric methods in
    algebra and number theory}, Progr. Math. 235, Birkh\"auser Boston,
  Boston, MA, 2005.

\bibitem{BCG} I.C. Bauer, F. Catanese \and F. Grunewald, The
  classification of surfaces with $p_{g} =q= 0$ isogenous to a product
  of curves, {\em Pure Appl. Math. Q. } {\bf 4}(2) (2008), 547--586.

\bibitem{magma} W. Bosma, J. Cannon \and C. Playoust, The Magma
  algebra system. I. The user language, {\em J. Symbolic Comput.}
  {\bf 24}(3-4) (1997), 235--265.

\bibitem{Bruck} R. H. Bruck, Quadratic extensions of cyclic planes, in {\em Proceedings Symposium in Applied Math.}, (1960) 15--44.

\bibitem{CMSZ1} D.I. Cartwright, A.M. Mantero, T. Steger \and
  A. Zappa, Groups acting simply transitively on the vertices of a
  building of type $\widetilde{A}_{2}$, I, {\em Geom. Dedicata} {\bf
    47}(2) (1993), 143--166.

\bibitem{CMSZ2} D.I. Cartwright, A.M. Mantero, T. Steger \and
  A. Zappa, Groups acting simply transitively on the vertices of a
  building of type $\widetilde{A}_{2}$, II, {\em Geom. Dedicata} {\bf
    47}(2) (1993), 167--223.

\bibitem{Cat00} F. Catanese, Fibred surfaces, varieties isogenous to a
  product and related moduli spaces, {\em Amer. J. Math.} 122 (1)
  (2000), 1--44.

\bibitem{EH} M. Edjvet \and J. Howie, Star graphs, projective planes
  and free subgroups in small cancellation groups, {\em Proc. London
    Math. Soc.} {\bf 57}(2) (1988), 301--328.

\bibitem{EV} M. Edjvet \and A. Vdovina, On the SQ-universality of
  groups with special presentations, {\em J. Group Theory} {\bf 13}(6)
  (2010), 923--931.

\bibitem{How} J. Howie, On the SQ-universality of $T(6)$-groups', {\em
    Forum Math.} {\bf 1}(3) (1989), 251--272.
    
\bibitem{HP} D. R. Hughes \and F. C. Piper, `Projective Planes', {\em
    Springer} New York (1973)
    
\bibitem{IV} S. Immervol and A. Vdovina, Partitions of projective
  planes and construction of polyhedra, Max-Planck-Institut f\:{u}r
  Mathematik, Bonn, Preprint Series 23 (2001).
 
\bibitem{LSV} A. Lubotzky, B. Samuels \and U. Vishne, `Explicit
  construction of Ramanujan complexes of type $\widetilde{A}_{d}$',
  {\em European J. Combin.} 26(6) (2005) 965--993.
 
\bibitem{LS} R.C. Lyndon \and P.E. Schupp, Combinatorial group theory,
  Classics in Mathematics, Reprint of the 1977 edition,
  Springer-Verlag, Berlin, 2001.
  
\bibitem{PV} N. Peyerimhoff \and A. Vdovina, Cayley graph
  expanders and groups of finite width, to appear in {\em J. Pure and
    Applied Algebra} (2011).

\bibitem{Ser} F. Serrano, Isotrivial fibred surfaces, {\em
    Ann. Mat. Pura Appl. (4)} {\bf 171} (1996), 63--81.
    
\end{thebibliography}
\end{document}